\documentclass[12pt,a4paper]{article}
\usepackage{latexsym}

\usepackage[leqno]{amsmath}
\usepackage{amssymb,amsfonts,amscd,amsthm,amsgen}
\usepackage[matrix,arrow,curve]{xy}
\usepackage{graphicx}

\numberwithin{equation}{section}

\newtheorem{thm}[equation]{Theorem}
\newtheorem{prop}[equation]{Proposition}
\newtheorem{cor}[equation]{Corollary}
\newtheorem{lem}[equation]{Lemma}
\newtheorem{specialcase}[equation]{Special case}

\theoremstyle{definition}
\newtheorem{definition}[equation]{Definition}
\newtheorem{example}[equation]{Example}
\newtheorem{problem}[equation]{Problem}

\newtheorem*{questionI}{Question I}
\newtheorem*{questionII}{Question II}
\newtheorem*{questionII'}{Question II$^\prime$}

\newtheorem{tble}[equation]{Table}

\theoremstyle{remark}
\newtheorem{remark}[equation]{Remark}

\newcommand{\point}{\operatorname{point}}

\newcommand{\incl}{\operatorname{incl}}

\newcommand{\id}{\operatorname{id}}

\newcommand{\im}{\operatorname{im}}

\newcommand{\scirc}{{\scriptstyle{\circ}}}

\begin{document}
\bibliographystyle{alpha}

\title{Coincidence free pairs of maps}

\author{Ulrich Koschorke}
\maketitle

\begin{abstract}
This paper centers around two basic problems of topological coincidence theory. First, try to measure (with help of Nielsen and minimum numbers) how far a given pair of ´maps is from being loose, i.e.\ from being homotopic to a pair of coincidence free maps. Secondly, describe the set of loose pairs of homotopy classes. We give a brief (and necessarily very incomplete) survey of some old and new advances concerning the first problem. Then we attack the second problem mainly in the setting of homotopy groups. This leads also to a very natural filtration of all homotopy sets. Explicit calculations are carried out for maps into spheres and projective spaces.
\end{abstract}

2000 {\it Mathematics Subject Classification.} Primary 55 M 20. Secondary 55 Q 40, 57 R 22.

{\it Key words and phrases.} Coincidence; Nielsen number; minimum number; configuration space; projective space, filtration.


%

\section{Introduction}                                                                          

Let \ $M$ \ be a closed smooth \ $m$--dimensional manifold.

In the first half of the 1920's S.\ Lefschetz established his celebrated results on fixed point theory which, in particular, yield the following.

\begin{thm} \                                                           
Let \ $f : M \to M$ \ be a map. If \ $L (f) \ne 0$, then every map \ $f'$ \ homotopic to \ $f$ \ has at least one fixed point \ $x \in M$ \ (i.e. $f' (x) = x)$.
\end{thm}

Here the Lefschetz number \ $L (f)$ \ can be defined to be the intersection number of the graph of \ $f$ \ with the diagonal \ $\Delta$ \ in \ $M \times M$.

The theorem of Lefschetz was a groundbreaking achievement. Still, it left several questions open.

\begin{questionI} \ What is the minimum number of fixed points of maps which are homotopic to \ $f$?

This is the principal problem of topological fixed point theory (cf.\ \cite{B}, p.\ 9).
\end{questionI}

\begin{questionII} \ What can we say about the set of homotopy classes \ $[f]$ \ which contain a fixed point free selfmap \ $f$ of \ $M$ \ (apart from the necessary condition \ $L (f) = 0$)?
\end{questionII}

In 1927 the Danish mathematician J.\ Nielsen achieved decisive progress concerning the first of these questions by decomposing the fixed point set of \ $f$ \ as follows. Two fixed points are called {\em (Nielsen) equivalent} \ if they can be joined by a path \ $\sigma$ \ in \ $M$ \ which is homotopic to \ $f \scirc \sigma$ \ by a homotopy which keeps the endpoints fixed. Each of the resulting Nielsen equivalence classes of fixed points contributes (trivially or nontrivially) to the Lefschetz number \ $L (f)$. Define \ $N (f)$ \ to be the number of {\em essential} \ Nielsen classes (i.e.\ of those classes whose contribution to \ $L (f)$ \ is not zero). This {\em Nielsen number} \ $N (f)$ \ is a lower bound for the {\bf m}inimum number of {\bf f}ixed points
\begin{equation}                                                             
MF (f) \ := \ \ \underset{f' \sim f}{\min} \ (\# \{ x \in M | f' (x) = x \} )
\end{equation}
and agrees with it except when \ $m = 2$ \ and the Euler number \ $\chi (M)$ \ of \ $M$ \ is strictly negative; however, in this exceptional case \ $MF (f) - N (f)$ \ can be arbitrarily large for suitable \ $f$ \ (compare \cite{N}, \cite{W}, \cite{Ji 1}, \cite{Ji 2}, \cite{Ji 3}, and \cite{Ke}; an excellent survey is given by R.\ Brown \cite{B}).

\smallskip
Fixed point questions allow a very natural and interesting generalization. Given a pair of maps \ $f_1, f_2 : M \to N$ \ between smooth connected manifolds (where \ $M$ \ is closed), let
\begin{equation}                                 
C (f_1, f_2) \ := \ \{ c \in M | f_1 (x) = f_2 (x)\}  \ \ \ \subset \ \ \ M
\end{equation}
denote its coincidence set. Here the relevant {\bf m}inimum number of {\bf c}oincidence points is
\begin{equation}                                  
MC (f_1, f_2) \ := \ \min \{ \# C (f'_1, f'_2) | f'_1 \sim f_1, f'_2 \sim f_2 \} \ .
\end{equation}

According to a result of R.B.S.\ Brooks \cite{Br} the same minimum number is achieved if we vary only one of the maps \ $f_1$ \ or \ $f_2$ \ by a homotopy; in particular, \ $MF (f) = MC (f, \text{identity map})$ \ (compare 1.2).

Unlike fixed point theory, coincidence theory allows the domain \ $M$ \ and the target \ $N$ \ to be different manifolds of arbitrary positive dimensions \ $m$ \ and \ $n$, resp. If \ $m > n$, then \ $MC (f_1, f_2)$ \ will be often infinite. Thus it makes sense to consider also the {\bf m}inimum number of {\bf c}oincidence {\bf c}omponents
\begin{equation}                                              
MCC (f_1, f_2) := \min \{ \# \pi_0 (C (f'_1, f'_2)) | f'_1 \sim f_1; f'_2 \sim f_2 \}
\end{equation}
which is always finite.

\begin{remark}                                                   
 (i) \ As one can see rather easily the values of the minimum numbers \ $MC$ \ and \ $MCC$ \ remain unchanged if we replace the base point free maps and homotopies in 1.4 and 1.5 by base point preserving ones (requiring e.g.\ that \ $f_1 (*) = *_1 \ne *_2 = f_2 (*)$ \ where \ $* \in M$ \ and \ $*_1, *_2 \in N$ \ are given base points; cf.\ e.g.\ the  appendix of \cite{K 6}).

(ii) \ Clearly \ $MC (f_1, f_2) = MCC (f_1, f_2) = 0$ \ whenever \ $m < n$ \  (use an approximation of \ $(f_1, f_2) : M \to N \times N$ \ which is transverse to the diagonal).
\end{remark}

\bigskip
In a series of recent papers (see, in particular, \cite{K 6} and \cite{K 7}) we studied the minimum numbers \ $MC$ \ and \ $MCC$ \ in arbitrary codimensions \ $m - n$. For this purpose we introduced a \lq\lq strong\rq\rq \ Nielsen number \ $N^\# (f_1, f_2)$ \ which generalizes Nielsen's original definition. Our approach is based on a careful analysis of the case when the pair \ $(f_1, f_2)$ \ is generic. Here the coincidence locus \ $C (f_1, f_2)$ \ is a smooth closed \ $(m - n)$--dimensional submanifold of \ $M$, equipped with

(i) \ \ a canonical description of its (nonstabilized) normal bundle; \ \ and

(ii) \ the map \ $\widetilde g$ \ from \ $C (f_1, f_2)$ \ to  the path space
\begin{equation}                                                
E (f_1, f_2) := \{ (x, \theta) \in M \times P (N) | \theta (0) = f_1 (x), \theta (1) = f_2 (x) \}
\end{equation}
defined by
\begin{equation*}
\widetilde g (x) = (x, \text{constant path at} \ f_1 (x) = f_2 (x)),
\end{equation*}
$x \in C (f_1, f_2)$, \ (where \ $P (N)$ \ denotes the space of all continuous paths \ $\theta : I \to N)$.

The space \ $E (f_1, f_2)$ \ has a very rich topology. Already its set \ $\pi_0 (E (f_1, f_2))$ \ of pathcomponents can be huge (it is  bijectively related to all wellstudied \lq\lq Reidemeister set\rq\rq; cf.\ \cite{K 3}, 2.1). The corresponding decomposition of \ $C (f_1, f_2)$ \ into the inverse images (under \ $\widetilde g$) of these pathcomponents generalizes the Nielsen decomposition of fixed point sets. But in higher codimensions \ $m - n > 0$ \ the map \ $\widetilde g$ \ into \ $E (f_1, f_2)$ \ can capture much further and deeper geometric information (which -- surprising often -- is related to various strong versions of Hopf invariants; see e.g.\ \cite{K 3}, 1.14, or \cite{K 6}, 7.6 and 7.11).

Details of the definition of \ $N^\# (f_1, f_2)$ \ can be found in \cite{K 6} where we proved also the following result.

\begin{thm}                                               
Let \ $f_1, f_2 : M^m \to N^n$ \ be {\rm (}continuous{\rm )} \ maps between smooth connected manifolds of the indicated dimensions, \ $M$ \ being closed. Then:
\begin{enumerate}
\item[\rm (i)] \ The Nielsen number \ $N^\# (f_1, f_2) = N^\# (f_2, f_1)$ \ is finite and depends only on the homotopy classes of \ $f_1$ \ and \ $f_2$.
\item[\rm (ii)] \ $0 \le N^\# (f_1, f_2) \le MCC (f_1, f_2) \le MC (f_1, f_2) \le \infty$; \ if \ $n \ne 2$ \ then also \ $MCC (f_1, f_2) \le \# \pi_0 (E (f_1, f_2))$; \ if \ $(m, n) \ne (2, 2)$ \ then \ $MC (f_1, f_2) \le \# \pi_0 (E (f_1, f_2))$ \ or \ $MC (f_1, f_2) = \infty$.
\item[\rm (iii)] \ If \ $M = N$ \ and \ $f_2 = \text{identity map}$, then \ $N^\# (f_1, f_2)$ \ coincides with the  Nielsen number of \ $f_1$ \ as defined in classical fixed point theory.
\end{enumerate}
\end{thm}

This allows us to compute our minimum numbers explicitly in various concrete cases.

\begin{example}  {\bf: \ spherical maps into spheres} \                               
(compare e.g. \cite{K 6}, 1.24). Consider maps \ $f_1, f_2 : S^m \to S^n$ \ where \ $m, n \ge 1$, \ and let \ $A$ \ denote the antipodal involution. Then
$$
N^{\#} (f_1, f_2) = MCC (f_1, f_2) =
\begin{cases}  0           & \text{if} \ f_1 \sim A {\scirc} f_2 \quad  ; \\
{\#} \pi_0 (E (f_1, f_2)) & \text{else} \qquad  .
\end{cases}
$$

If \ $f_1 \not\sim A {\scirc} f_2$ \ then \ $\# \pi_0 (E (f_1, f_2))$ \ equals \ $1$ \ (and \ $|d^0 (f_1) - d^0 (f_2)|$, \ resp.)\ according as \ $n \ne 1$ \ (or \ $m = n = 1$, resp.; here \ $d^0 (f_i)$ \ denotes the usual degree).

Clearly \ $MC (f_1, f_2) \le 1$ \ whenever \ $[f_1] - [A {\scirc} f_2]$ \ lies in \ $E (\pi_{m -1} (S^{n -1}))$, the image of the Freudenthal suspension. On the other hand, it is wellknown that \ $MC (f_1, f_2)$ \ is infinite if  \ $[f_1] - [A {\scirc} f_2] \not\in E (\pi_{m -1} (S^{n -1}))$ \ and \ $(m, n) \ne (1, 1)$.
\end{example}

\begin{example} {\bf : \ spherical maps into projective spaces} (cf.\ \cite{K 7}, 1.17).
Let \ $\mathbb K \ = \ \mathbb R, \ \mathbb C$ \ or \ $\mathbb H$ \ denote the field of real, complex or quaternionic numbers, and let \ $d = 1, 2$ \ or \ $4$ \ be its real dimension. Let \ $\mathbb K P (n')$ \ and \ $V_{n' + 1, 2} (\mathbb K)$, resp., denote the corresponding space of lines and orthonormal 2--frames, resp., in \ $\mathbb K^{n' +1}$. The real dimension of \ $N = \mathbb K P (n')$ \ is \ $n := d \cdot n'$. Consider the diagram
\begin{equation}
\begin{CD} \cdots \to \pi_m (V_{n' + 1, 2} (\mathbb K)) @>{p_{\mathbb K *}}>> \pi_m (S^{n + d - 1}) @>{\partial_{\mathbb K}}>> \pi_{m -1} (S^{n -1}) \to \cdots \\
@.                     @VV{p_*}V                           @VV{E}V \\
@.                       \pi_m (\mathbb K P (n'))      @.        \pi_m (S^n)
\end{CD}
\end{equation}
determined by the canonical fibrations \ $p$ \ and \ $p_{\mathbb K}$; \ $E$ \ denotes the Freudenthal suspension homomorphism.

We want to determine the minimum numbers of all pairs of maps \ $f_1, f_2 : S^m \to \mathbb K P (n'), \ m, n' \ge 1$. In view of example 1.9 and remark 1.6 (ii) we need not consider the cases where \ $n' = 1$ \ (or \ $m = 1$).
\end{example}

\begin{lem}                                         
Assume \ $m, n' \ge 2$. Then \ $p_*$ \ (cf.\ 1.11) is injective and
$$
\pi_m (\mathbb K P (n')) = p_* (\pi_m (S^{n + d - 1})) \oplus \pi_m^c (\mathbb K P (n'))
$$
where \ $\pi_m^c (\mathbb K P (n')) := \incl_* (\pi_m (\mathbb K P (n') - \{ *\}))$ \ and \ $\incl$ \ denotes the inclusion of \ $\mathbb K P (n')$, punctured at some point \ $*$. Hence, given \ $[f_i] \in \pi_m (\mathbb K P (n'))$, \ there is a unique homotopy class \ $[\widetilde f_i] \in \pi_m (S^{n + d - 1})$ \ such that \ $p_* ([\widetilde f_i]) - [f_i] \in \pi_m^c (\mathbb K P (n')), \ i = 1, 2$. (Since \ $\pi_m^c (\mathbb K P (n')) \cong \pi_{m -1} (S^{d -1})$, \ we may assume that \ $\widetilde f_i$ \ is a genuine lifting of \ $f_i$ \ when \ $\mathbb K = \mathbb R$ \ or when \ $m > 2$ \ and \ $\mathbb K = \mathbb C$).
\end{lem}

We see this by comparing the exact homotopy sequences of the fibrations
\begin{equation}                                  
p \ : \ S^{n + d - 1} \ \ \longrightarrow \ \ \mathbb K P (n')
\end{equation}
and \ $p| : S^{n -1} \to \mathbb K P (n' - 1) \ (\underset{\subset}{\sim} \mathbb K P (n') - \{ *\} )$.

\begin{thm} \                                                 
Assume \ $m, n' \ge 2$. Each pair of homotopy classes \ $[f_1], [f_2] \in \pi_m (\mathbb K P (n'))$ \ satisfies precisely one of the seven conditions which are listed in table 1.15, together with the corresponding Nielsen and minimum numbers. {\rm(}Here  we use the language of lemma {\rm 1.12} and define also \ $[f'_i] := [p \scirc \widetilde f_i] \in \pi_m (\mathbb K P (n')), \ i = 1, 2$; \ moreover \ $A$ \ denotes the antipodal map on \ $S^{n + d - 1})$.
\end{thm}

\renewcommand{\arraystretch}{1.5}
\begin{tabular}{l|c|c|c|} \hline
{Condition } & ${\scriptstyle N^\# (f_1, f_2)}$ & ${\scriptstyle MCC (f_1, f_2)}$ & ${\scriptstyle MC (f_1, f_2)}$ \\ \hline\hline
1) \ ${\textstyle f'_1 \, \sim \, f'_2,\ \ [\widetilde f_2] \ \in \ \ker \partial_{\mathbb K}}$  & 0 & 0 & 0 \\ \hline
2) \ ${\textstyle f'_1 \, \sim \, f'_2,\ \ [\widetilde f_2] \ \in \ \ker E \scirc \partial_{\mathbb K} -  \ker \partial_{\mathbb K}}$ & 0 & 1 & 1 \\ \hline
3) \ ${\textstyle \mathbb K \, = \, \mathbb R, \ f'_1 \sim f'_2, \ \ \widetilde f_2 \not\sim A \scirc \widetilde f_2}$ & 1 & 1 & 1 \\ \hline
4) \ ${\textstyle  \mathbb K  =  \mathbb R,} \ {\scriptstyle f'_1 \not\sim f'_2, \ \ [\widetilde f_1] - [\widetilde f_2]}  \in  {\textstyle E}\, {\scriptstyle (\pi_{m -1} (S^{n -1}))} $ & 2 & 2 & 2 \\ \hline
5) \ ${\textstyle \mathbb K \, = \, \mathbb R,  \ [\widetilde f_1] - [\widetilde f_2] \ \not\in \ E (\pi_{m -1} (S^{n -1}))}$ & 2 & 2 & $\infty$ \\ \hline
6) \ ${\textstyle \mathbb K \, = \, \mathbb C \, \text{or} \, \mathbb H, \ [\widetilde f_1] \, = \, [\widetilde f_2] \ \not\in \ \ker E \scirc \partial_{\mathbb K}}$ & 1 & 1 & 1 \\ \hline
7) \ ${\textstyle \mathbb K \, = \, \mathbb C \, \text{or} \, \mathbb H, \ [\widetilde f_1] \ne [\widetilde f_2]}$ & 1 & 1 & $\infty$ \\ \hline
\end{tabular}

\begin{tble} \   
Nielsen and minimum coincidence numbers of all pairs of maps \ $f_1, f_2 : S^m \to \mathbb K P (n'), \ m, n' \ge 2$: \ replace each (possibly base point free) homotopy class \ $[f_i]$ \ by a base point preserving representative and read off the values of \ $N^\#$ \ and \ $M(C)C$. (Here \ $f'_1 \sim f'_2$ \ means that \ $f'_1, f'_2$ \ are homotopic in the base point free sense. For proofs see \cite{K 7}).

\smallskip
This concludes our brief (and necessarily rather incomplete) survey of some of the developments triggered by our initial Question I. \hfill $\square$
\end{tble}

In this paper we start investigating a natural generalization of Question II.

\begin{definition}                              
Let \ $M$ \ and \ $N$ \ be smooth connected manifolds, \ $M$ \ being closed. A pair of maps \ $f_1, f_2 : M \to N$ \ is called {\em loose} \ if it is homotopic to a coincidence free pair; in other words, if \ $MC (f_1, f_2) = 0$ \ or, equivalently \ $MCC (f_1, f_2) = 0$.

\smallskip
It makes no difference whether we use base point free or base point preserving homotopies in this definition (provided \ $f_1 (*) \ne f_2 (*)$ \ when \ $*$ \ is a given base point of \ $M$; \ cf.\ 1.6).
\end{definition}

\begin{questionII'}
What can we say about the set of homotopy classes of loose pairs?
\end{questionII'}

We will concentrate on the case \ $M = S^m, \ m \ge 1$. Let \ $* \in S^m$ \ and \ $*_1 \ne *_2 \in N$ \ be given base points.

Consider the subgroups
$$
\pi^c_m (N, *_i) \ \ \ \subset \ \ \ \pi^{(2)}_m (N, *_i) \ \ \ \subset \ \ \ \pi_m (N, *_i), \ \ \ \ \ i = 1,2,
$$
where
\medskip
\begin{equation}                            
\begin{aligned}
\pi^c_m &(N, *_i) := \{ [f] \in \pi_m (N, *_i) | (f, *_{i \pm 1}) \ \text{is loose} \}  =  \incl_* (\pi_m (N - \{ *_{i \pm 1} \},  \ *_i)) \\
\hfill  \text{and} \\
\pi^{(2)}_m &(N, *_i) := \{ [f] \in \pi_m (N, *_i) | \ \exists [\overline f] \in \pi_m (N, *_{i \pm 1}) \ \text{s. t.} \ (f, \overline f) \ \text{is loose} \}.
\end{aligned}
\end{equation}
Here \ $*_i$ \ denotes also the constant map with the indicated value, and \ $\incl$ \ stands for the obvious inclusion. (Compare also remark 3.7 below).

\begin{thm}                            
For \ $m \ge 1$ \ there is a welldefined group isomorphism
$$
c \ \ : \ \ \pi^{(2)}_m (N, *_1) \big/ \pi^c_m (N, *_1) \ \ \ \longrightarrow \ \ \ \pi^{(2)} (N, *_2) \big/ \pi^c_m (N, *_2)
$$
which takes the coset \ $[[f]]$ \ of \ $[f] \in \pi^{(2)}_m (N, *_1)$ \ to the coset of any element \ $[\overline f] \in \pi^{(2)}_m (N, *_2)$ \ such that \ $(f, \overline f)$ \ is loose.

A pair \ $([f_1], [f_2]) \in \pi_m (N, *_1) \times \pi_m (N, *_2)$ \ is loose if and only if \ $[f_i] \in \pi^{(2)}_m (N, *_i), \ i = 1, 2$, \ and \ $c ([[f_1]]) = [[f_2]]$.

In particular, if \ $\pi_m^c (N, *_i) = \pi_m (N, *_i)$ \ then all pairs of maps \ $f_1, f_2 : S^m \to N$ \ {\rm (}base point preserving or not{\rm)} are loose; this is the case e.g.\ when \ $N$ \ is not compact or when \ $m < n$.
\end{thm}

\begin{specialcase}                         
If \ $N$ \ allows a fixed point free selfmap \ $A: \ N \to N$ \ such that \ $A (*_1) = *_2$ \ then \ $\pi^{(2)}_m (N, *_i) = \pi_m (N, *_i)$ \ for all \ $m \ge 1$, and \ $c$ \ is induced by \ $A$ {\rm (}i.e. $c([[f]]) = [[A \scirc f]])$.
\end{specialcase}

Thus (the nontriviality of) \ $\pi_m (N, *_1) / \pi^{(2)}_m (N, *_1)$ \ is an obstruction to the existence of such a fixed point free selfmap. On the other hand, such selfmaps occur e.g.\ on the total space of every nontrivial covering map.

\begin{example}  ${\bf (N = S^n):}$                             
A pair of maps \ $f_1, f_2 : S^m \to S^n$ \ (base point preserving or not) is loose if and only if \ $f_1 \sim A \scirc f_2$ \ where \ $A$ \ denotes the antipodal map. Indeed, \ $\pi^c_m (S^n) = \{ 0\}$ \ and
$$
\begin{CD}
c \ = \ A_* \ : \ \pi_m (S^n,  *_1 ) \ @>{\cong}>> \ \pi_m (S^n,  A (*_1))
\end{CD}
$$
(compare also \cite{DG}, 2.10).
\end{example}

\begin{cor}                      
If \ $N$ \ allows a nowhere vanishing vector field (e.g.\ if \ $N$ \ is odd); then
$$
\pi^{(2)}_m (N, *_i) \ = \ \pi_m (N, *_i) \ \ \ \text{for all} \ m \ge 1, \ \ i = 1, 2 \ ,
$$
and \ $c$ is induced by a map \ $A$ {\rm(} as in {\rm 1.19)} \ which is homotopic to the identity.
\end{cor}

Indeed, the flow of the vector field yields the required fixed point free map \ $A$.

Theorem 1.18 and corollary 1.21 suggest that Question II$^\prime$ \ may be most interesting when \ $N$ \ is a closed even--dimensional manifold.

\begin{example} {\bf (Projective spaces)}.            
Consider the case \ $N =  \mathbb K P (n'), \ \mathbb K \ = \ \mathbb R, \ \mathbb C$ \ or \ $\mathbb H, \ \ m, n' \ge 2$, \ as in 1.10 (and use the language of lemma 1.12).  Then a pair of maps \ $f_1, f_2 : S^m \to \mathbb K P (n')$ \ is loose precisely if the corresponding pair \ $(p \scirc \widetilde f_1, p \scirc \widetilde f_2)$ \ is loose or, equivalently, if the maps \ $p \scirc \widetilde f_1, p \scirc \widetilde f_2$ \ are homotopic and \ $\widetilde f_i : S^m \to S^{n + d - 1}$ \ can be lifted to the Stiefel manifold \ $V_{n' + 1, 2} (\mathbb K), \ i = 1 $  or  $2$ \ (compare 1.11). Thus here the isomorphism \ $c$ \ (cf.\ 1.18) is induced by a selfmap \ $A$ \ of \ $\mathbb K P (n')$ \ which is homotopic to the identity map but which can be fixed point free only when \ $\mathbb K \ = \ \mathbb R$ \ and \ $n$ \ is odd, i.e.\ when the Lefschetz number \ $L (A) = \chi (\mathbb K P (n'))$ \ vanishes. \hfill    $\square$
\end{example}

\begin{problem}                                    
Is the group isomorphism \ $c$ \ in theorem 1.18 always induced by a selfmap of \ $N$?
\end{problem}

\medskip

It seems to be very desirable to determine \ $\pi^c_* (N), \ \pi^{(2)}_* (N), \ c$ \ and hence the sets of loose pairs of homotopy classes (cf.\ theorem 1.18) for many more concretely given closed sample manifolds, e.g.\ for Stiefel manifolds and Grassmannians. Here is a partial result in this direction.

\begin{example}                             
For every {\em even} \ integer \ $r \ge 4$, all pairs of maps \ $f_1, f_2 : S^m \to G_{r, 2} (\mathbb R)$ \ into the Grassmann manifold of 2--planes in \ $\mathbb R^r$ \ are loose.
\end{example}

\medskip

Details of the proof of theorem 1.18 and of its consequences will be given in section 2.

\medskip
Throughout our discussion a central role is played by the set of homotopy classes of those maps which occur in loose pairs (i.e.\ which are  {\em not coincidence producing} \  in the terminology of Brown and Schirmer, cf.\ \cite{BS}). For arbitrary topological spaces \ $X$ \ and \ $Y$ \ this set turns out to be the first interesting term of a very natural descending filtration of the full homotopy set \ $[X, Y]$. In section 3 we study this filtration and determine it e.g.\ for the homotopy groups of spheres and projective spaces.

\section{Loose pairs}                      

Throughout this paper manifolds are required to be Hausdorff spaces having a countable basis and no boundary.

\begin{proof}[Proof of theorem 1.18] \
If the pairs \ $(f, \overline f), (f, \overset {=}{f})$ \ and \ $(\hat f, *_2)$ \ are loose, then so are \ $(*_1 \sim f \cdot f^{- 1}, \ \overline f \cdot \overset{=}{f}^{- 1})$ \ and \ $(f \cdot \hat f \ , \overline f \cdot *_2 \sim \overline f)$. Thus the coset \ $[[\overline f]] = [[\overset{=}{f}]]$ \ is determined by \ $[f] \in \pi^{(2)}_m (N, *_1)$; \ it does not depend on the choice of the class \ $[\overline f]$ \ $($which  makes \ $(f, \overline f)$ \ loose$)$ and not even on the choice of \ $[f]$ \ within its coset.

If the pairs \ $(f, \overline f)$ \ and \ $(f', \overline f')$ \ are loose, then so is also \ $(f \cdot f', \overline f \cdot \overline f')$. Hence the bijection \ $c$ \ which interchanges the roles of \ $f_1$ \ and \ $f_2$ \ in a loose pair \ $(f_1, f_2)$ \ is compatible with the group structure.

If \ $N$ \ is not compact or if \ $m < n$, then every map \ $f : S^m \to N$ \ can be deformed into the complement of a given point in \ $N$ via a suitable isotopy or via transverse approximation.
\end{proof}

Next we turn to the situation of example 1.22. Given arbitrary (not necessarily base point preserving) maps \ $f_i : S^m \to \mathbb K P (n')$, \ put \ $*_i := f_i (*)$ \ and choose \ $\widetilde f_i$ \ as in lemma 1.12; then the summand \ $[p \scirc \widetilde f_i] - [f_i]$ \ plays no role in looseness questions, \ $i = 1, 2$. If the pair \ $(p \scirc \widetilde f_1, p \scirc \widetilde f_2)$ \ is coincidence free, then the unit vectors \ $\widetilde f_1 (x), \widetilde f_2 (x)$ \ in \ $\mathbb K^{n' + 1}$ \ are linearly independent for all \ $x \in S^m$; \ suitable rotations yield both a homotopy \ $\widetilde f_1 \sim \widetilde f_2$ \ and liftings to the Stiefel manifold of {\em orthonormal} \ 2--frames in \ $\mathbb K^{n' +1}$.

This argument shows also that the isomorphism
\begin{equation}                                         
\pi_m (\mathbb K P (n')) \ / \ \pi_m^c (\mathbb K P (n')) \ \ \ \cong \ \ \ \pi_m (S^{n + d - 1})
\end{equation}
induced by \ $p$ \ (c.f.\ 1.12 and 1.13) makes \ $\pi^{(2)}_m (\mathbb K P (n')) / \pi^c_m (\mathbb K P (n'))$ \ correspond to \ $\im (p_{\mathbb K *}) = \ker \partial_{\mathbb K}$ \ in diagram 1.11. This yields an alternative proof of the calculations in table 1.15 as far as condition 1) is concerned. Furthermore there are easy examples (e.g.\ when \ $\mathbb K = \mathbb R$ \ and \ $m = n' \equiv 0 (2)$, cf.\ 3.13) where \ $\partial_{\mathbb K}$ \ and hence \ $\pi_m (\mathbb K P (n')) / \pi^{(2)}_m (\mathbb K P (n'))$ \ is nontrivial  so that every selfmap of \ $\mathbb K P (n')$ \ must have a fixed point (compare 1.19). Of course such questions can be settled more systematically by the Lefschetz fixed point theorem. \hfill $\square$

\smallskip
Finally we prove the statement in example 1.24. In view of 1.6 (ii) we may assume that \ $m \ge 4$. According to theorem 1.18 our claim is established once we see that
$$
\incl_* : \pi_m (G_{r, 2} (\mathbb R) - \{ \point\} ) \ \longrightarrow \ \pi_m (G_{r, 2} (\mathbb R))
$$
is surjective. But this follows from

\begin{lem} \                                    
For all \ $m \ge 3$ \ and for all even integers \ $r = 2r' \ge 4$ \ we have the isomorphism
$$
\begin{CD}
e_* + u_* : \pi_m (\mathbb R P (r - 2)) \oplus \pi_m (\mathbb C P (r' - 1)) \ @>{\cong}>>  \ \pi_m (G_{r, 2} (\mathbb R))
\end{CD}
$$
where \ $e (\lambda) = \lambda \oplus \mathbb R (0, \dots, 0, 1), \ \lambda \in \mathbb R P (r - 2)$, and \ $u$ \ assigns the underlying real plane to any complex line. (Note that both embeddings \ $e$ \ and \ $u$ \ have codimensions \ $r - 2 > 0$).
\end{lem}

\begin{proof}
Scalar multiplication with the complex number \ $i$ \ on \ $\mathbb C^{r'} \cong \mathbb R^r$ determines a section \ $s$ \ of the fibration \ $S^{r -2} \subset V_{r, 2} (\mathbb R) \to S^{r -1}$.

Thus the exact homotopy sequence splits and yields the isomorphism
$$
\widetilde e_* + s_* \ : \ \pi_m (S^{r -2}) \oplus \pi_m (S^{r -1}) \ \longrightarrow \pi_m (V_{r, 2} (\mathbb R)) \ .
$$
This implies our claim since the fiber \ $O (2)$ \ of the projection \ $V_{r, 2} (\mathbb R) \ \to \ G_{r, 2} (\mathbb R)$ \ is aspherical.
\end{proof}

\begin{problem}                         
What is known about the groups \ $\pi_*^c$ \ and \ $\pi_*^{(2)}$ \ of arbitrary Grassmannians \ $G_{r, k} (\mathbb K), \ \ r - 1 > k > 1$ \ ?
\end{problem}

Note that in the special case \ $r = 2k$ \ there is the fixed point free involution \ $\perp$ \ on \ $G_{2k, k} (\mathbb K), \ \ \mathbb K \ = \ \mathbb R, \ \mathbb C$ \ or \ $\mathbb H$ \  (take orthogonal complements). Thus \ $\pi_*^{(2)} (G_{2 k, k} (\mathbb K))$ \ is the full homotopy group.

\section{A filtration of homotopy sets}             

In this section we extend the definition of the group \ $\pi_*^{(2)} (N)$ \ (formed by those maps which occur in loose pairs) and obtain a very natural infinite descending filtration of arbitrary homotopy sets \ $[M, N]$.

Given any topological space \ $N$, consider the commuting diagram

\begin{equation}                 
  \xymatrix{
    N = \widetilde{C}_1(N) \ar[rd]_{\id=p_1} &
    \ar[l] \widetilde{C}_2(N) \ar[d]^{p_2} &
    \ar[l] \cdots &
    \ar[l] \widetilde{C}_q(N) \ar[dll]^{p_q} & \ar[l] \cdots \\
    & N
  }
\end{equation}
of  configuration spaces
$$
\widetilde C_q (N) \ = \ \{ (y_1, \dots, y_q) \in N^q | \ y_i \ne
y_j \ \text{for} \ 1 \le i \ne j \le q \ \}
$$
$q \ge 1$, and of projections which drop the last component(s) of an (ordered)
configuration \ $(y_1, \dots, y_q)$. For any topological space \
$M$ \ this leads to the filtration

\begin{equation}            
[M, N] = [M, N]^{(1)} \supset [M, N]^{(2)} \supset \dots \supset
[M, N]^{(q)} := p_{q *} ([M, \widetilde C_q (N)])  \supset \dots
\end{equation}

Thus \ $[M, N]^{(q)}$ \ consists of the homotopy classes of those
maps which fit into a $q$-tuple \ $f_1, \dots, f_q : M \to N$ \ of
maps without any (pairwise) coincidences.

Next, given a base point \ $* \in M$ \ and an infinite sequence \
$(*_1, *_2, \dots)$ \ of pairwise distinct points in \ $N$, \
equip \ $\widetilde C_q (N)$ \ with the base point \ $(*_1, *_2,
\dots, *_q)$ \ and consider also the base point preserving version
of the filtration 3.2.

\begin{example}[$N = S^n, \ n \ge 1$] \             
For every point \ $y
\in S^n$ \ use the stereographic projection \ $\sigma_y$ \ from \
$S^n - \{y\}$ \ to the orthogonal complement of the line \ $\mathbb R
y$ in $\mathbb R^{n +1}$ \ (i.e.\ to the tangent space \ $T_y (S^n))$
\ and obtain the following fiber preserving homeomorphisms and
homotopy equivalences over \ $S^n$:
$$
\begin{aligned} \widetilde C_2 (S^n) \ \ & = \ \ S^n \times S^n - \Delta
\ \ \cong \ \ TS^n \ \ \qquad \qquad
 \ \ \text{and} \\
\widetilde C_3 (S^n) \ \ & \sim \ \ TS^n - \ \text{zero section} \
\ \overset{\supset}{\sim}  \ \ V_{n + 1,2} \ \ .
\end{aligned}
$$
Here e.g.\ the vectors \ $0$ \ and \ $v \ne 0$ \ in \ $T_y (S^n)$
\ correspond to the configurations \ $(y, -y) \in \widetilde C_2
(S^n)$ \ and \ $(y, -y, \sigma_y^{- 1} (v)) \in \widetilde C_3
(S^n)$.

When \ $q \ge 3$ \ the projection \ $\widetilde C_q (S^n) \to
\widetilde C_3 (S^n)$ \ has a section which corresponds to the map
\ $v \to (v, 2v, \dots, (q - 2) v)$.

Thus both in the base point free and in the base point preserving
setting we have for every topological space \ $M$ and \ $q \ge 3$
$$
[M, S^n] = [M, S^n]^{(2)} \supset [M, S^n]^{(3)} = p_{\mathbb R *} ([M, V_{n +
1, 2} (\mathbb R)]) =  [M, S^n]^{(q)} .
$$
As in 1.11 \ $p_{\mathbb R} : V_{n + 1,2} (\mathbb R) \to S^n$ \ denotes the standard projection
from the Stiefel manifold \ $ST (S^n)$ \ of unit tangent vectors.
If \ $n$ \ is odd it has a section and \ $[M, S^n]^{(q)} = [M,
S^n]$ \ for \ {\it all} \ $q \ge 1$. However, if \ $n$ \ is even
and e.g.\ $M = S^n$, then $[M, S^n]^{(2)} \ne [M, S^n]^{(3)}$.
\end{example}

\begin{prop} \                
 Both in the base point free and base
point preserving setting we have
$$
[M, N]^{(q)} \ = \ [M, N] \ \ \ \text{for all} \ q \ge 1
$$
if at least one of the following condition holds:

\parindent0pt {\rm(i)} \ \
$M$ is compact, but \ $N$ \ is not -- in the base point preserving
setting we assume also that \ $N$ \ is a connected topological
manifold; \ \ or

\parindent0pt
 {\rm(ii)} \ \ $N$ is a smooth manifold which allows a
nowhere vanishing vector field; \ \ or

\parindent0pt
 {\rm(iii)} \ $M$ \ and \ $N$ \ are smooth manifolds such
that \ $\dim M < \dim N$.
\end{prop}
\medskip

\begin{proof}

\parindent0pt (i) \  \ For every map \ $f : M \to N$ \
the complement \ $N - f (M)$ \ contains infinitely many points.
Thus for \ $q \ge 2$ \ there exist (e.g.\ constant) maps \ $f_2,
\dots, f_q$ \ which, together with \ $f_1 = f$, \ define the
required map into \ $\widetilde C_q (N)$. If \ $N$ \ is a
connected topological manifold we may assume that \ $f_i (*) =
*_i, \ i =1, \dots, q$, \ e.g.\ after suitable isotopies.

\parindent0pt
 (ii) \  \ We use the resulting flow \ $\varphi$ \ and a
suitable function \ $\varepsilon : N \longrightarrow (0, \infty)$
\ to define the pairwise coincidence-free selfmaps $ \id = A_1,
A_2, \dots, A_q$ \ of \ $N$ \ by
$$
A_i (x) \ = \ \varphi (x, (i - 1) \cdot \varepsilon (x) / q) , \ \
\ x \in N .
$$

Then \ $(f, A_2 \scirc f, \dots, A_q \scirc f)$ \ is a lifting of
\ $f$ \ to \ $\widetilde C_q (N)$; suitable modifications (e.g.\
by finger moves) allow us to make it base point preserving.

\parindent0pt
(iii) \ After a transverse approximation \ $f$ \ maps \
$M$ \ into \ $N - \{ *_2, \dots, *_q\}$.
\end{proof}

We will be mainly interested in the case \ $M = S^m, \ m \ge 1$
(studied in the base point preserving setting). We obtain the
nested sequence of subgroups
\begin{equation}                       
\pi_m (N, *_1) = \pi^{(1)}_m (N, *_1) \supset \pi^{(2)}_m (N, *_1)
\supset \dots \supset \pi^{(q)}_m (N, *_1) \supset \dots
\end{equation}
defined by
\begin{equation}                       
\pi^{(q)}_m (N, *_1) \ := \ p_{q *} (\pi_m (\widetilde C_q (N),
(*_1, \dots, *_q))), \ \ \ q\ge 1 .
\end{equation}

For \ $q = 2$ \ this agrees with the definition in 1.17 since \
$p_2$ \ is the first projection on \ $\widetilde C_2 (N) = N
\times N - \Delta$.

\begin{remark} \              
Assume that \ $N$ is a topological
manifold of dimension \ $n \ge 1$.

Then all the arrows in diagram 3.1 are projections of locally
trivial fibrations (compare \cite{FN}). It follows from the homotopy lifting property  that, given a loose
pair \ $(f_1, f_2)$, \ only one of the maps \ $f_i$, \ say \
$f_2$, \ has to be deformed to \ $f'_2$ \ so that \ $(f_1, f'_2)$
\ is coincidence free (compare \cite{Br}). In particular for all \
$m \ge 1$
\medskip
\begin{equation} \             
\pi^c_m (N, *_1) \ = \ \incl_* (\pi_m (N - \{ *_2\}, *_1)) \
\subset \ \bigcap_{q \ge 1} \pi^{(q)}_m (N, *_1) \ =: \ \pi^{(\infty)}_m (N, *_1)
\end{equation}
(compare 1.17 and 3.5; here \ $\incl$ \ denotes the inclusion of \ $N
- \{*_2\} \cong N$--small ball around \ $*_2$).

Moreover we have the exact sequence
\medskip
\begin{equation}                      
\begin{CD}
\dots \to \pi_m (N \times N - \Delta, (*_1, *_2)) @>{ \ p_{2 *} \ }>>
\pi_m (N, *_1) @>{ \ \delta_m \ }>> \pi_{m -1} (N - \{*_1\}, *_2)
\end{CD}
\end{equation}
where \ $p_2$ \ denotes the projection \ $(y_1, y_2)
\to y_1$. Thus
\begin{equation}                   
\pi^{(2)}_m (N, *_1) \ = \ \ker \delta_m \ .
\end{equation}
In view of proposition 3.4 we will be particularly interested in the case where \ $N$ \ is a closed connected manifold of even dimension \ $n \le m$.
\end{remark}

\noindent
{\bf Example 3.11:} \ $N = \mathbb K P (n')$ \ where \ $\mathbb K = \mathbb R, \ \mathbb C$ \ or \ $\mathbb H$          
(compare 1.10). \ In view of 3.3 and 3.4 (iii) we need to consider only the case \ $m, n' \ge 2$. Thus we can (and will) use the terminology of 1.11 and 1.12.

\stepcounter{equation}
\begin{prop} \                                  
Assume \ $n' \ge 2$. Then for all \ $q \ge 2$ \ {\rm(} and \ $m \ge 1${\rm )}
$$
\pi^{(q)}_m (\mathbb K P (n')) \ = \ \pi^{(2)}_m (\mathbb K P (n')) \ = \ p_* (\ker \partial_{\mathbb K}) \ \oplus \  \pi^c_m (\mathbb K P (n')) \ ;
$$
the analogous result holds for base point free homotopy classes of arbitrary maps \ $f : S^m \to \mathbb K P (n')$.

Moreover let \ $M$ \ be any paracompact space. If \ $\mathbb K \ = \ \mathbb R$ \ and \ $H^1 (M; \mathbb Z_2) = 0$, or if \ $\mathbb K \ = \ \mathbb C$ \ and \ $H^2 (M; \mathbb Z) = 0$, then for all \ $q \ge 2$
$$
[M, \mathbb K P (n')]^{(q)} \ = \ [M, \mathbb K P (n')]^{(2)} \ ;
$$
this set coincides with the full homotopy set \ $[M, \mathbb K P (n')]$ \ if in addition \ $n'$ \ is odd.
\end{prop}

\begin{proof} \
In each case we need to consider only maps \ $f$ \ which allow a lifting \ $\widetilde f : M \to \ S^{n + d - 1}$, i.e. $f = p \scirc \widetilde f$ \ (compare 1.13). If \ $M = S^m$ \ this follows from 3.8; otherwise use characteristic classes to see that the pullback of the canonical line bundle over \ $\mathbb K P (n')$ \ is trivial.

Given liftings \ $\widetilde f, \overset{\simeq}{f}$ \ such that the pair \ $(p \scirc \widetilde f, p \scirc \overset{\simeq}{f})$ \ is coincidence free, \ $\widetilde f (x), \overset {\simeq}{f} (x)$ \ are linearly independent unit vectors in \ $\mathbb K^{n' + 1}$ \ for all \ $x \in M$.
Thus \ $p \scirc \widetilde f$ \ is the starting term of an (arbitrarily long) sequence of pairwise coincidence free maps \ $f_i : M \to \mathbb K P (n')$ \ defined by
$$
f_i (x) = \mathbb K (\widetilde f (x) + (i - 1) \overset{\simeq}{f} (x)), \ \ x \in M, \ i \ge 1 \ .
$$
We conclude that \ $[M, \mathbb K P (n')]^{(2)} \subset [M, \mathbb K P (n')]^{(q)}$.

If \ $n'$ \ is odd and \ $\mathbb K = \mathbb R$ \ or \ $\mathbb C$, then \ $p_{\mathbb K}$ \ (cf.\ 1.11) allows a section (via the complex or quaternionic scalar multiplication on \ $\mathbb K^{n' + 1}$) and every map \ $f = p \scirc \widetilde f$ \ occurs in a coincidence free pair as above.
\end{proof}

In contrast, when \ $n'$ \ is even then \ $[M, \mathbb K P (n')]^{(2)}$ \ often turns out to be strictly smaller than the full homotopy set \ $[M, \mathbb K P (n')]$ \ (or, in the terminology of Brown and Schirmer, there are {\em coincidence producing} \ maps from \ $M$ \ to \ $\mathbb K P (n')$, cf.\ \cite{BS}). Let us illustrate this for \ $\mathbb K = \mathbb R$.

\begin{lem} \                                    
For all \ $m, n > 1$ \ the diagram
$$
\xymatrix@C=5pc{
 \pi_m(S^n) \ar[r]^{\partial_{\mathbb R}} \ar[rd]_(0.4){(1+(-1)^n)E^\infty \; \; \; \; \;} & \pi_{m-1}(S^{n-1}) \ar[d]^{E^\infty}\\
 & \pi^S_{m - n}
}
$$
commutes up to a fixed sign. {\rm (}Here \ $E^\infty$ \ denotes stable
suspension.{\rm )}

In particular, in the stable dimension range \ $m < 2n - 2$ \ {\rm (}where both arrows labelled \ $E^\infty$ \ are isomorphisms{\rm )} we have
$$
\ker \partial_{\mathbb R} = \{ z \in \pi_m (S^n) | \ (1 + (- 1)^n) \cdot z = 0\} \ .
$$
\end{lem}

\begin{proof} \ Given \ $[f] \in \pi_m (S^n)$, \ the Freudenthal suspension of \ $\partial_{\mathbb R} ([f])$ \ equals the selfintersection invariant \ $\pm \, \underline\omega (f, f)$ \ (cf.\ \cite{K 7}, 5.6 and 5.7). In turn
$$
E^\infty (\underline\omega (f, f)) \ = \ \omega (f, f) \ = \ \pm \chi (S^n) \cdot E^\infty ([f])
$$
in \ $\Omega_{m - n} (S^m; \varphi) \cong \pi^S_{m -n}$ \ (cf.\ \cite{K 2}, 1.4, 1.6, and 2.2; here we use the canonical stable trivializations  of the tangent bundle \ $TS^n$ \ and of the virtual coefficient bundle \ $\varphi = f^* (TS^n) - TS^m)$.
\end{proof}

\begin{example}                           
$$
\begin{aligned}
 \pi_9^c (\mathbb R P (6)) &= 0 \ \  \ \subset \ \
\pi_9^{(2)} (\mathbb R P (6)) \cong \mathbb Z_2 \ \ \ \subset \ \ \pi_9
(\mathbb R P (6)) \ \cong \
\mathbb Z_{24} \\
\text{and} \qquad \ \ \ \ & \\
\pi_{17}^c (\mathbb R P (10)) &= 0 \ \ \ \subset \ \ \pi_{17}^{(2)}
(\mathbb R P (10)) \cong \mathbb Z_2 \ \ \subset \ \ \pi_{17} (\mathbb R P
(10)) \cong \mathbb Z_{240}
\end{aligned}
$$

This follows from our results 1.12, 3.12, 3.13, and Toda's tables
\cite{T}.
\end{example}

\begin{remark}                         
Assume that \ $N$ \ is a topological manifold. For \ $i = 1, 2$ \ consider the fiber inclusion and the projection
$$
\begin{CD}
(N - \{ *_i\}, \ *_{i \pm 1}) \ @>{\subset}>{\textstyle j_i}> \ (\widetilde C_2 (N), (*_1, *_2)) @>>{\textstyle p_{\scriptstyle 2, i}}> (N, *_i)
\end{CD}
$$
of the locally trivial fibration defined by \ $p_{2, i} (y_1, y_2) = y_i$; its exact homotopy sequence (cf.\ 3.9) yields the isomorphism
$$
\begin{CD}
p_{2, i*} : \pi_m (\widetilde C_2 (N), (*_1, *_2)) / (\im j_{1 *} + \im j_{2 *}) @>{\cong}>>  \pi^{(2)}_m (N, *_i) / \pi^c_m (N, *_i) .
\end{CD}
$$
Then the composite \ $p_{2, 2*} \scirc p^{- 1}_{2, 1*}$ \ equals the group isomorphism \ $c$ \ (cf.\ theorem 1.18) which is so central to our study of loose pairs of maps.
\end{remark}

\parindent0pt
Universit\"at Siegen \\ Emmy Noether Campus, Walter-Flex-Str. 3,  D-57068 Siegen, Germany \\
E-mail: koschorke@mathematik.uni-siegen.de\\
url: http://www.math.uni-siegen.de/topology

\end{document}